\documentclass[11pt]{amsart}
\usepackage{pstricks}
\usepackage{latexsym}
\usepackage{amssymb}

 \textwidth=\textwidth \textheight=\textheight
\theoremstyle{plain}
\newtheorem{thm}{Theorem}
\newtheorem{lem}[thm]{Lemma}
\newtheorem{prop}[thm]{Proposition}
\newtheorem{lemma}[thm]{Lemma}
\newtheorem{corollary}[thm]{Corollary}

\newtheorem{conj}[thm]{Conjecture}

\theoremstyle{definition}
\newtheorem{definition}[thm]{Definition}

\newtheorem{conjdef}[thm]{Conjectural Definition}

\newtheorem{example}[thm]{Example}

\newcommand{\set}[1]{ \left\{ #1 \right\} }

\def\P{{\mathcal P}}
\def\u{{\mathbf u}}
\def\U{{\mathcal U}}

\def\C{\mathbb C}
\def\ll{\lambda}
\def\tll{\tilde{\lambda}}

\def\p{{\mathcal P}}

\def\zz{{\mathbb Z}}
\def\Z{{\mathbb Z}}

\def\uqsln{U_q(\widehat{ \mathfrak{
sl}}_n)}

\def\g{{\mathcal G}}
\def\f{{\mathbf F}}

\def\nn{{\mathbb N}}

\def\reading{\mathrm{reading}}

\def\spin{\mathrm{spin}}

\newcommand{\tbt}[4]{\begin{array}{cc} j_#1 & j_#2 \\ j_#3 &
j_#4\end{array}}

\newcommand{\be}[1]{\begin{equation} \label{#1}}
\newcommand{\ee}{\end{equation}}

\newcommand{\br}[1]{\langle #1\rangle}
\newcommand{\brac}[1]{\left(#1\right)}
\newcommand{\h}[2]{h_{#1}^{(#2)}}

\psset{unit=1pt} \psset{arrowsize=4pt 1} \psset{linewidth=1pt}
\newrgbcolor{mygreen}{.2 .7 .2}
\newrgbcolor{myred}{.7 .3 .2}
\newgray{mygray}{.90}
\countdef\x=23 \countdef\y=24 \countdef\z=25 \countdef\t=26

\def\abox(#1,#2)#3{
\x=#1 \y=#2 \multiply\x by 16 \multiply\y by 16 \z=\x \t=\y
\advance\z by 16 \advance\t by 16
\psline(\x,\y)(\x,\t)(\z,\t)(\z,\y)(\x,\y) \advance\x by 8
\advance\y by 8 \rput(\x,\y){{\bf #3}}}

\def\hdom(#1,#2)#3{
\x=#1 \y=#2 \multiply\x by 16 \multiply\y by 16 \z=\x \t=\y
\advance\z by 32 \advance\t by 16
\psline(\x,\y)(\x,\t)(\z,\t)(\z,\y)(\x,\y) \advance\x by 16
\advance\y by 8 \rput(\x,\y){{\bf #3}}}

\def\vdom(#1,#2)#3{
\x=#1 \y=#2 \multiply\x by 16 \multiply\y by 16 \z=\x \t=\y
\advance\z by 16 \advance\t by 32
\psline(\x,\y)(\x,\t)(\z,\t)(\z,\y)(\x,\y) \advance\x by 8
\advance\y by 16 \rput(\x,\y){{\bf #3}}}

\def\rec(#1,#2,#3,#4){
\psline(#1,#2)(#3,#2)(#3,#4)(#1,#4)(#1,#2) }

\begin{document}
\title{Ribbon Schur Operators}
\author{Thomas Lam}
\address{Department of Mathematics,
         M.I.T., Cambridge, MA 02139}

\email{thomasl@math.mit.edu}
\date{September 2004}
\begin{abstract}
A new combinatorial approach to the ribbon tableaux generating
functions and $q$-Littlewood Richardson coefficients of Lascoux,
Leclerc and Thibon \cite{LLT} is suggested.  We define operators
which add ribbons to partitions and following Fomin and Greene
\cite{FG} study non-commutative symmetric functions in these
operators. This allows us to give combinatorial interpretations
for some (skew) $q$-Littlewood Richardson coefficients whose
non-negativity appears not to be known.  Our set up also leads to
a new proof of the action of the Heisenberg algebra on the Fock
space of $\uqsln$ due to Kashiwara, Miwa and Stern \cite{KMS}.
\end{abstract}
\maketitle
\section{Introduction}
The aim of this paper is to provide a new combinatorial framework
to study ribbon tableaux and the ribbon tableaux generating
functions $\g_{\ll/\mu}(X;q)$ of Lascoux, Leclerc and Thibon
\cite{LLT}.  The functions $\g_{\ll/\mu}(X;q)$ are to ribbon
tableaux what skew Schur functions are to usual (skew) Young
tableaux.

Lascoux, Leclerc and Thibon developed the theory of ribbon
tableaux generating functions in connection with the Fock space
representation $\f$ of $\uqsln$.  Though the definition of
$\g_{\ll/\mu}(X;q)$ is purely combinatorial, Lascoux, Leclerc and
Thibon showed that $\g_{\ll/\mu}(X;q)$ were symmetric functions
using a theorem of \cite{KMS}, which describes an action of the
Heisenberg algebra on $\f$.  The theory was further developed in
\cite{LT} by showing that the $q$-Littlewood Richardson
coefficients $c^\nu_\ll(q)$ given by
\[
\g_{\ll}(X;q) = \sum_\nu c^\nu_\ll(q) s_\nu(X)
\]
were related to the canonical basis of $\f$.  More recently,
Haglund et. al. \cite{HHLRU} and Haglund \cite{Hag} have
conjectured mysterious relationships between the ribbon tableaux
generating functions on the one hand and Macdonald polynomials and
diagonal harmonics on the other hand.

We will study ribbon tableaux generating functions in an
elementary fashion, proving their symmetry, the result of
Kashiwara, Miwa and Stern \cite{KMS} concerning the action of the
Heisenberg algebra on $\f$, and also developing an approach to the
non-negativity of the polynomials $c^\nu_{\ll/\mu}(q)$ (known to
be true for $\mu = \emptyset$, \cite{LT}).  Our approach follows
Fomin and Greene's work on non-commutative Schur functions
\cite{FG}.  We define \emph{ribbon Schur operators} $u_i$ which
add a ribbon to a partition with head along diagonal $i$, if
possible, and multiplying by a power of $q$ according to the
\emph{spin} of the ribbon.  These operators are related to the
LLT-ribbon tableaux generating functions in the same way that the
operators which add boxes along specified diagonals are related to
Schur functions (see \cite{FG}).  Following \cite{FG}, we study
non-commutative symmetric functions, in particular non-commutative
Schur functions $s_\nu(\u)$ in the operators $u_i$.  We show
directly that $s_\nu(\u)$ is a non-negative sum of monomials for
the cases when $\nu$ is a hook shape or of the form $\nu = (s,2)$.
This in turn proves the non-negativity of the polynomials
$c^\nu_{\ll/\mu}(q)$ for any skew shape $\ll/\mu$, and $\nu$ of
the above form.

In Section \ref{sec:ribSchur} we begin by defining semistandard
ribbon tableaux and the ribbon Schur operators $u_i$.  In Section
\ref{sec:algebra}, we describe all the relations in the algebra
generated by the ribbon Schur operators.  In Section
\ref{sec:symmetry}, we show that the non-commutative
``homogeneous'' symmetric functions $h_k(\u)$ in the operators
$u_i$ commute, and define the non-commutative Schur functions
$s_\nu(\u)$.  In Section \ref{sec:cauchy}, we prove a Cauchy style
identity for the operators $u_i$ and their adjoints $d_i$.  This
also proves the action of the Heisenberg algebra on $\f$ of
\cite{KMS}, and as shown in \cite{Lam}, it implies ribbon Pieri
and Cauchy formulae for ribbon tableaux generating functions.  In
Section \ref{sec:ribbonfunctions}, we relate the non-commutative
Schur functions $s_\nu(\u)$ to the $q$-Littlewood Richardson
coefficients $c^\nu_{\ll/\mu}(q)$.  In Section \ref{sec:schur}, we
give positive combinatorial descriptions of $s_\nu(\u)$ for the
cases where $\nu$ is a hook shape or a shape of the form $(a,2)$
or $(2,2,1^a)$.  We conjecture that $s_\nu(\u)$ always has such a
positive description.  Finally, in Section \ref{sec:final}, we
give some concluding remarks concerning ribbon Schur operators.

\bigskip
{\bf Acknowledgements.}  This project is part of my Ph.D.
Thesis written under the guidance of Richard Stanley.  I am
grateful for all his advice and support over the last couple of
years.  I would also like to thank Sergey Fomin for interesting
discussions and Marc Van Leeuwen for comments on parts of this
work.

\section{Ribbon Schur operators}
\label{sec:ribSchur} We will follow mainly \cite{Mac} and
\cite{EC2} for notation related to partitions and symmetric
functions.

Let $\ll$ be a partition.  We will always draw our partitions in
the English notation.  The \emph{diagonal} or \emph{content} of a
square $(i,j) \in \ll$ is equal to $j-i$. Let $n$ be a fixed
positive integer. A $n$-ribbon is a connected skew shape of size
$n$ containing no $2 \times 2$ square. The \emph{head} of a ribbon
is the box at its top-right.  The \emph{spin} $\spin(R)$ of a
ribbon $R$ is equal to the number of rows in the ribbon, minus 1.

Let $K$ denote the field $\C(q)$.  Let $\p$ denote the set of
partitions and $\f$ denote a vector space over $K$ spanned by a
countable basis $\set{\ll \mid \ll \in \p}$ indexed by partitions.
Define linear operators $u_i^{(n)}: \f \rightarrow \f$ for $i \in
\Z$ which we call \emph{ribbon Schur operators} by:
\[
u_i^{(n)}: \ll \longmapsto \begin{cases} q^{\spin(\mu/\ll)} \mu &
\mbox{if
$\mu/\ll$ is a $n$-ribbon with head lying on the $i$-th diagonal,} \\
    0 &  \mbox{otherwise.} \end{cases}
\]
We will usually suppress the integer $n$ in the notation, even
though $u_i$ depends on $n$.  If we need to emphasize this
dependence, we write $u_i^{(n)}$.  We say that a partition $\ll$
has an \emph{$i$-addable ribbon} if a $n$-ribbon can be added to
$\ll$ with head on the $i$-th diagonal.  Similarly, $\ll$ has an
\emph{$i$-removable ribbon} if a $n$-ribbon can be removed from
$\ll$ with head on the $i$-th diagonal.

\begin{figure}[ht]
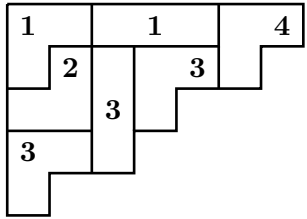

\pspicture(0,20)(120,80)

\psline(0,90)(112,90)(112,74)(96,74)(96,58)(64,58)(64,42)(48,42)(48,26)(16,26)(16,10)(0,10)(0,90)
\psline(0,58)(16,58)(16,74)(80,74)(80,90) \psline(80,74)(80,58)
\psline(32,90)(32,26) \psline(48,74)(48,42) \psline(0,42)(32,42)
\rput(8,34){{\bf 3}} \rput(24,66){{\bf 2}} \rput(40,50){{\bf 3}}
\rput(72,66){{\bf 3}} \rput(56,82){{\bf 1}} \rput(8,82){{\bf 1}}
\rput(104,82){{\bf 4}}

\endpspicture
\caption{A semistandard $3$-ribbon tableau with shape
$(7,6,4,3,1)$, weight $(2,1,3,1)$ and spin 7.}
\label{fig:ribbontableau}
\end{figure}

A skew shape is called a \emph{horizontal ribbon strip} if it can
be tiled with ribbons so that the heads of every ribbon touch the
top of the skew shape.  If such a tiling exists, it is necessarily
unique and will be called the \emph{natural} tiling.
Alternatively, $\ll/\mu$ is a horizontal ribbon strip if there
exists $i_1 < i_2 < \cdots < i_k$ so that
$u_{i_k}u_{i_{k-1}}\cdots u_{i_2}u_{i_1}. \mu = q^a \ll$ for some
integer $a$.  A \emph{semistandard ribbon tableau} $T$ of shape
$\ll/\mu$ is a chain $(\mu = \ll^0 \subset \ll^1 \subset \cdots
\subset \ll^r = \ll)$ of partitions where $\ll^i/\ll^{i-1}$ is a
horizontal ribbon strip for each $i \in [1,r]$.  We usually
represent $T$ by tiling each $\ll^i/\ll^{i-1}$ with ribbons, in
the natural way, and filling each ribbon with the integer $i$ (see
Figure \ref{fig:ribbontableau}). The \emph{spin} $\spin(T)$ of a
tableau $T$ is the sum of the spins of its composite ribbons.

Littlewood's \emph{$n$-quotient bijection} is a bijection between
partitions $\ll$ and its \emph{$n$-quotient}
$(\ll^{0},\ll^{(1)},\ldots,\ll^{(n-1)}) \in \P^n$ together with
$n$-core $\tll$.  This bijection is described in \cite{Mac, LLT}.
We will describe it in terms of the operators $u_i$. The $n$-core
$\tll$ is obtained from $\ll$ by removing $n$-ribbons until it is
no longer possible to do so.  Each $n$-core $\delta$ is determined
(not uniquely) by an integer vector $(s_0,s_1,\ldots,s_{n-1}) \in
\Z^n$ of \emph{offsets} which satisfies $s_i \neq s_j \mod n$ for
$i \neq j$.  The offsets interact with the $n$-quotient in the
following way. Suppose $\ll$ and $\mu$ are two partitions with
$\tll = \delta = \tilde{\mu}$ and so that $\ll^{(i)} = \mu^{(i)}$
for all $i \neq j$ for some $j \in \set{0,1,\ldots,n-1}$.  If
$\mu^{(j)} = u^{(1)}_k \ll^{(j)}$ then $\mu = u^{(n)}_{nk+s_j}
\ll$.  In other words, the operators $\set{u_i^{(n)} \mid i = s_j
\ \mod n}$ act on the $j$-th partition of the $n$-quotient
$\ll^{(j)}$ by adding squares. This together with the fact that
the $n$-quotient $(\emptyset,\emptyset, \ldots,\emptyset)$ of an
$n$-core is empty determines the bijection.


\medskip
Lascoux, Leclerc and Thibon \cite{LLT} defined weight generating
functions $\g_{\ll/\mu}(X;q)$ (from hereon called \emph{ribbon
functions}) for ribbon tableaux:
\[
\g_{\ll/\mu}(X;q) = \sum_T q^{\spin(T)} x^T
\]
where the sum is over all semistandard ribbon tableaux of shape
$\ll/\mu$.  They proved that these functions were symmetric
functions and defined the $q$-Littlewood Richardson coefficients
$c^\nu_{\ll/\mu}(q)$ by
\[
\g_{\ll/\mu}(X;q) = \sum_\nu c^\nu_{\ll/\mu}(q) s_\nu(X).
\]
In \cite{LT}, it was shown that $c^\nu_{\ll/\mu}(q)$ was a
non-negative polynomial in $q$ for $\mu = \emptyset$, and the
proof there can be generalised to $\mu$ being any $n$-core (see
\cite{HHLRU}). However, it appears that the non-negativity is not
known for general skew shapes $\ll/\mu$.

\section{The algebra of ribbon Schur operators}
\label{sec:algebra} Fix an integer $n \geq 1$ throughout. Let $\U
= \U_n \subset \mathrm{End}_K(\f)$ denote the algebra generated by
the operators $\set{u_i^{(n)}}$ over $K$.

\begin{prop}
\label{prop:commute} The operators $u_i$ satisfy the following
commutation relations:
\begin{align}
\label{eq:rel1}
u_i u_j &= u_j u_i &\mbox{for $|i - j| \geq n+1$,} \\
\label{eq:rel2}
u_i^2 & = 0& \mbox{for $i \in \Z$,}\\
\label{eq:rel3}
u_{i+n}u_iu_{i+n} &= 0 & \mbox{for $i \in \Z$,}\\
\label{eq:rel4}
u_i u_{i+n} u_i &= 0& \mbox{for $i \in \Z$,}  \\
\label{eq:rel5} u_i u_j &= q^2 u_j u_i &\mbox{for $n > i - j > 0$}.
\end{align}
Furthermore, these relations generate all the relations that the
$u_i$ satisfy.  The algebra generated by $\set{u_i}_{i=1}^{i=kn}$
has dimension $(C_k)^n$ where $C_k = \frac{1}{2k+1}{2k \choose k}$
is the $k$-th Catalan number.
\end{prop}

\begin{figure}
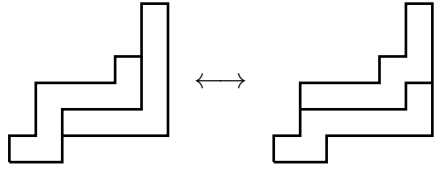

\pspicture(0,0)(160,60)
\psline(0,0)(0,10)(10,10)(10,30)(40,30)(40,40)(50,40)(50,60)(60,60)(60,10)(20,10)(20,0)(0,0)
\psline(20,10)(20,20)(50,20)(50,40)
\rput(80,30){$\longleftrightarrow$}
\psline(100,0)(100,10)(110,10)(110,30)(140,30)(140,40)(150,40)(150,60)(160,60)(160,10)(120,10)(120,0)(100,0)
\psline(110,20)(150,20)(150,30)(160,30)
\endpspicture
 \caption{Calculating relation (\ref{eq:rel5})
of Proposition \ref{prop:commute}.} \label{fig:spinswap}
\end{figure}

\begin{proof} Relations (\ref{eq:rel1}-\ref{eq:rel4}) follow from
the description of ribbon tableaux in terms of the $n$-quotient
and the usual relations for the operators $u_i^{(1)}$ (see for
example \cite{FG}).  Relation (\ref{eq:rel5}) is a quick
calculation (see Figure \ref{fig:spinswap}), and also follows from
the inversion statistics of \cite{SSW,HHLRU} which give the spin
in terms of the $n$-quotient.

Now we show that these are the only relations.  The usual Young
tableau case with $n=1$ was shown by Billey, Jockush and Stanley
in \cite{BJS}. However, when $q = 1$, we are reduced to a direct
product of $n$ copies of this action as described earlier.  The
operators $\set{u_i \mid i = k \mod n}$ act on the $k$-th tableau
of the $n$-quotient independently.

Let $f = a_u \underline{u} + a_v \underline{v} + \cdots \in \U$
and suppose $f$ acts identically as 0 on $\f$.  First suppose that
some monomial $\underline{u}$ acts identically as $0$.  Then by
the earlier remarks, the subword of $\underline{u}$ consisting
only of $\set{u_i \mid i=k \mod n}$ must act identically as 0 for
some $k$.  Using relation (\ref{eq:rel5}), we see that we can
deduce $\underline{u} = 0$, using the result of \cite{BJS}.

Now suppose that a monomial $\underline{v}$ does not act
identically as 0, so that $\underline{v} \cdot \mu = q^t \ll$ for
some $t \in \Z$ and $\mu, \ll \in \P$.  Collect all other
monomials $\underline{v}'$ such that $\underline{v}' \cdot \mu =
q^{b(\underline{v}')+ t} \ll$ for some $b(\underline{v}') \in \Z$.
By Lemma \ref{lem:commute} below, $\underline{v}' =
q^{b(\underline{v}')}  \underline{v}$.  Since $f \cdot \mu = 0$ we
must have $\sum_{\underline{v}'} a_{v'} \underline{v}' = 0$ and by
Lemma \ref{lem:commute} this can be deduced from the relations
(\ref{eq:rel1}-\ref{eq:rel5}).  This shows that we can deduce $f =
0$ from the relations.

For the last statement of the theorem, relation (\ref{eq:rel5})
reduces the statement to the case $n=1$.  When $n=1$, we think of
the $u_i$ as the Coxeter generators $s_i$ of $S_{k+1}$.  A basis
of the algebra generated by $\{u_i^{(1)}\}_{i=1}^k$ is given by
picking a reduced decomposition for each $321$-avoiding
permutation -- these are exactly the permutations with no
occurrences of $s_i s_{i+1} s_i$ in any reduced decomposition. It
is well known that the number of these permutations is equal to a
Catalan number.
\end{proof}
\begin{lem}
\label{lem:commute} Suppose $\underline{u} =
u_{i_k}u_{i_{k-1}}\cdots u_{i_1}$ and $\underline{v} =
u_{j_l}u_{j_{l-1}}\cdots u_{j_1}$.  If $\underline{u} \cdot \mu =
q^t \underline{v} \cdot \mu \neq 0$ for some $t \in \Z$ and $\mu
\in \P$ then $\underline{u} = q^t \underline{v}$ as operators on
$\f$, and this can be deduced from the relations of Proposition
\ref{prop:commute}.
\end{lem}
\begin{proof}
Using relation (\ref{eq:rel5}) and the $n$-quotient bijection, we
can reduce the claim to the case $n=1$ which we now assume.  So
suppose $\underline{u} \cdot\mu = \underline{v}\cdot \mu = \ll$.
Then the multiset of indices $\set{i_1,i_2,\ldots,i_k}$ and
$\set{j_1,j_2,\ldots,j_l}$ are identical, since these are the
diagonals of $\ll/\mu$.  In particular we have $k = l$.

We need only show that using relations
(\ref{eq:rel1}-\ref{eq:rel4}) we can reduce to the case where
$u_{i_1} = u_{j_1}$, and the result will follow by induction on
$k$.  Let $a = \min\set{b \mid j_b = i}$ where we set $i = i_1$.
We can move $u_{j_a}$ to the right of $\underline{v}$ unless for
some $c < a$, we have $j_c = i\pm 1$.  But $\mu$ has an
$i$-addable corner, and so $\nu = u_{j_{c-1}}u_{j_{c-2}}\cdots
u_{j_1} \cdot \mu$ also has an $i$-addable corner.  This implies
that $u_{i \pm 1} \cdot \nu = 0$ so no such $c$ can exist.
\end{proof}

\section{Non-commutative homogeneous symmetric functions in ribbon Schur
operators} \label{sec:symmetry} Let $h_k(\u) = \sum_{i_1 < i_2 <
\cdots < i_k} u_{i_k}\cdots u_{i_1}$ be the ``homogeneous''
symmetric functions in the operators $u_i$ (the name makes sense
since $u_i^2 = 0$). Since the $u_i$ do not commute, the ordering
of the variables is important in the definition.  The action of
$h_k(\u)$ on $\f$ is well defined though it does not lie in $\U$.
By the remarks in Section \ref{sec:ribSchur}, the operator
$h_k(\u)$ adds a horizontal ribbon strip of size $k$ to a
partition. Alternatively, we may write
\[
\prod_{i = \infty}^{i = -\infty} \left( 1 + xu_i \right) =
\sum_{i=0}^{\infty} x^k h_k(\u)
\]
where if $i < j$ then $(1+xu_i)$ appears to the right of
$(1+xu_j)$ in the product.

The following proposition was shown in \cite{LLT} using
representation theoretic results in \cite{KMS}.  Our new proof
imitates \cite{FG,Fom}.

\begin{thm}
\label{thm:hcommute} The elements $\set{h_k(\u)}_{k=1}^{\infty}$
commute and generate an algebra isomorphic to the algebra of
symmetric functions.
\end{thm}

\begin{proof}
Given any fixed partition only finitely many $u_i$ do not
annihilate it.  Since we are adding only a finite number of
ribbons, it suffices to prove that $h_k(u_a,u_{a+1},\ldots,u_b)$
commute for every two integers $b > a$.  First suppose that $ n
\geq b -a$.  We may assume without loss of generality that $a = 1$
and $b = n+1$.  We expand both
$h_k(u_a,u_{a+1},\ldots,u_b)h_l(u_a,u_{a+1},\ldots,u_b)$ and
$h_l(u_a,u_{a+1},\ldots,u_b)h_k(u_a,u_{a+1},\ldots,u_b)$ and
collect monomials with the same set of indices $I = \set{i_1 < i_2
< \cdots < i_{k+l}}$.  By Proposition \ref{prop:commute}, any
operator $u_i$ can occur at most once in any such monomial.
Suppose the collection of indices $I$ does not contain both 1 and
$n+1$. Then by Proposition \ref{prop:commute} we may reorder any
such monomial $\underline{u} = u_{j_1}u_{j_2}\cdots u_{j_{k+l}}$.
into the form $q^t u_{i_1}u_{i_2} \cdots u_{i_{k+l}} = q^t u_I$.
The integer $t$ is given by twice the number of inversions in the
word $j_1j_2\cdots j_{k+l}$. That the coefficient of $u_I$ is the
same in $h_k(u_a,u_{a+1},\ldots,u_b)h_l(u_a,u_{a+1},\ldots,u_b)$
and $h_l(u_a,u_{a+1},\ldots,u_b)h_k(u_a,u_{a+1},\ldots,u_b)$ is
equivalent to the following generating function identity for
permutations (alternatively, symmetry of a Gaussian polynomial).

Let $D_{m,k}$ be the set of permutations of $S_m$ with a single
ascent at the $k$-th position and let $d_{m,k}(q) = \sum_{w \in
D_{m,k}}q^{{\rm inv}(w)}$ where ${\rm inv}(w)$ denotes the number
of inversions in $w$. Then we need the identity
\[
d_{k+l,k}(q) = d_{k+l,l}(q).
\]
This identity follows immediately from the involution on
permutations $w = w_1w_2\cdots w_m \in S_m$ given by $w \mapsto v$
where $v_i = m+1 -w_{m+1-i}$.

When $I$ contains both $1$ and $n+1$ then we must split further
into cases depending on the locations of these two indices: (a)
$u_{n+1} \cdots u_1 \cdots$; (b) $\cdots u_{1} u_{n+1} \cdots$;
(c) $\cdots u_{n+1} \cdots u_1$; (d) $u_{n+1} \cdots u_1$.  We
pair case (a) of
$h_k(u_a,u_{a+1},\ldots,u_b)h_l(u_a,u_{a+1},\ldots,u_b)$ with case
(c) of $h_l(u_a,u_{a+1},\ldots,u_b)h_k(u_a,u_{a+1},\ldots,u_b)$
and vice versa; and also cases (b) and (d) with itself. After this
pairing, and using relation (\ref{eq:rel5}) of Proposition
\ref{prop:commute}, the argument goes as before.  For example, in
cases (a) and (c), we move $u_{n+1}$ to the front and $u_1$ to the
end.

Now we consider $h_k(u_a,\ldots,u_b)$ for $b-a > n$. Let
$E_{b,a}(x) = ( 1 + xu_b)(1+xu_{b-1}) \cdots (1+xu_a)$. Note that
$E_{b,a}(x)^{-1} = (1-xu_a)(1-x u_{a+1}) \cdots (1-xu_b)$ is a
valid element of $\U[x]$. The commuting of the
$h_k(u_a,u_{a+1},\ldots,u_b)$ is equivalent to the following
identity:
\[
E_{b,a}(x) E_{b,a}(y) = E_{b,a}(y) E_{b,a}(x)
\]
as power series in $x$ and $y$ with coefficients in $\U$ which we
assume to be known for all $(b,a)$ satisfying $b - a < l $ for
some $l > n$. Now let $b = a+ l$.  In the following we use the
fact that $u_a$ and $u_b$ commute.
\begin{align*}
&E_{b,a}(x) E_{b,a}(y) \\
&= E_{b,a+1}(x)(1+xu_a)(1+yu_b) E_{b-1,a}(y)\\
&=E_{b,a+1}(y)E_{b,a+1}(x)\brac{E_{b,a+1}(y)}^{-1}(1+yu_b)(1+xu_a)\brac{E_{b-1,a}(x)}^{-1}E_{b-1,a}(y)E_{b-1,a}(x)\\
&=E_{b,a+1}(y)E_{b,a+1}(x)\brac{E_{b-1,a+1}(x)E_{b-1,a+1}(y)}^{-1}E_{b-1,a}(y)E_{b-1,a}(x)\\
&=E_{b,a+1}(y)(1+xu_b)(1+yu_a)E_{b-1,a}(x)\\
&=E_{b,a}(y)E_{b,a}(x).\\
\end{align*}
This proves the inductive step and thus also that the $h_k(\u)$
commute.  To see that they are algebraically independent, we may
restrict our attention to an infinite subset of the operators
$\set{u_i}$ which all mutually commute, in which case the
$h_k(\u)$ are exactly the classical elementary symmetric functions
in those variables.
\end{proof}

Theorem \ref{thm:hcommute} allows us to make the following
definition, following \cite{FG}.

\begin{definition}
The \emph{non-commutative (skew) Schur functions}
$s_{\ll/\mu}(\u)$ are given by the Jacobi-Trudi formula:
\[
s_{\ll/\mu}(\u) = \det\left(h_{\ll_i - i +
j-\ll_j}(\u)\right)_{i,j=1}^{l(\ll)}.
\]
Similarly we may define the non-commutative symmetric function
$f(\u)$ for any symmetric function $f$.
\end{definition}

It is not clear at the moment how to write $s_\ll(\u)$, or even
just $e_k(\u) = s_{1^k}(\u)$, in terms of monomials in the $u_i$
like in the definition of $h_k(\u)$. One can check, for example,
that $p_2(\u)$ cannot be written as a non-negative sum of
monomials.


\section{The Cauchy identity for ribbon Schur operators}
\label{sec:cauchy}

The vector space $\f$ comes with a natural inner product
$\br{.,.}$ such that $\br{\ll,\mu} = \delta_{\ll \mu}$.  Let $d_i$
denote the adjoint operators to the $u_i$ with respect to this
inner product. They are given by
\[
d_i : \ll \longmapsto \begin{cases} q^{\spin(\ll/\mu)}\mu &
\mbox{if $\ll/\mu$ is a $n$-ribbon with head lying on the $i$-th
diagonal,} \\ 0 & \mbox{otherwise}.
\end{cases}
\]

The following lemma is a straightforward computation.
\begin{lemma}
\label{lem:udcommute} Let $i \neq j$ be integers.  Then $u_i d_j =
d_j u_i$.
\end{lemma}

Define $U(x)$ and $D(x)$ by
\[
U(x) = \cdots (1+xu_{2})(1+xu_1)(1+xu_0)(1+xu_{-1}) \cdots
\]
and
\[
D(x) = \cdots (1+xd_{-2})(1+xd_{-1})(1+xd_0)(1+xd_{1}) \cdots.
\]
So $U(x) = \sum_k x^k h_k(\u)$ and we similarly define
$h_k^\perp(\u)$ by $D(x) = \sum_k x^k h_k^\perp(\u)$.  The
operator $h_k^\perp(\u)$ acts by removing a horizontal ribbon
strip of length $k$ from a partition.

The main result of this section is the following identity.

\begin{thm}\label{thm:cauchy}
The following Cauchy Identity holds:
\begin{equation}
\label{eq:cauchy} U(x)D(y) \prod_{i=0}^{n-1} \frac{1}{1-q^{2i}xy}
= D(y)U(x).
\end{equation}
\end{thm}
A combinatorial proof of this identity was given by Marc van
Leeuwen \cite{vL} via an explicit shape datum for a
Schensted-correspondence.  Our proof is suggested by ideas in
\cite{Fom} but considerably different to the techniques there.
The following Corollary is immediate after equating coefficients
of $x^ay^b$ in Theorem \ref{thm:cauchy}.
\begin{corollary}
\label{cor:hhperp} Let $a,b \geq 1$ and $m = min(a,b)$.  Then
\[
h_b^\perp(\u) h_a(\u) = \sum_{i=0}^m h_i(1,q^2,\ldots,q^{2(n-1)})
h_{a-i}(\u)h_{b-i}^\perp(\u).
\]
\end{corollary}

Combining Theorem \ref{thm:cauchy} with Theorem \ref{thm:hcommute}
we obtain the following corollary, first proved in \cite{KMS} (the
connection with ribbon tableaux was first shown in \cite{LLT}):
\begin{corollary}
\label{cor:KMS} The operators $B_{-k} = p_k(\u)$ and $B_k =
p_k^\perp(\u)$ generate an action of the Heisenberg algebra on
$\f$ with commutation rule
\begin{equation}
\label{eq:hei} [B_k,B_l] =  k \cdot \frac{1 -
q^{2n|k|}}{1-q^{2|k|}} \cdot \delta_{k,-l}.
\end{equation}
\end{corollary}
\begin{proof}
When $k$ and $l$ have the same sign, the commutation relation
follows from Theorem \ref{thm:hcommute}.  For the other case, we
first write (see \cite{Mac,EC2})
\[
h_a(\u) = \sum_{\ll \vdash a} z_\ll^{-1}p_\ll(\u),
\]
where $z_\ll = 1^{m_1(\ll)}2^{m_2(\ll)}\cdots m_1(\ll)! m_2(\ll)!
\cdots$ and $m_i(\ll)$ denotes the number of parts of $\ll$ equal
to $i$.  We first show that (\ref{eq:hei}) implies Corollary
\ref{cor:hhperp}.  Thus we need to show that (\ref{eq:hei})
implies
\begin{multline*}
\brac{\sum_{\ll \vdash b}z_\ll^{-1}p_\ll^\perp(\u)}\brac{
\sum_{\ll \vdash a} z_\ll^{-1}p_\ll(\u)}  \\ =  \sum_{i=0}^m
h_i(1,q^2\ldots,q^{2(n-1)}) \brac{\sum_{\ll \vdash
 a-i}z_\ll^{-1} p_\ll(\u)}\brac{ \sum_{\ll \vdash b-i}
 z_\ll^{-1}p_{\ll}^\perp(\u)}.
\end{multline*}
Note that  $p_k(1,q,\ldots,q^{2(n-1)}) = \frac{1 -
q^{2n|k|}}{1-q^{2|k|}}$.  Let $\mu$ and $\nu$ be partitions such
that $m = |\nu|=|\mu|$.  One checks that the coefficient of
$p_\mu(\u) p_\nu^\perp(\u)$ on the right hand side is equal to
$z_\nu^{-1}z_\mu^{-1}\sum_{\ll \vdash m}
z_\ll^{-1}p_\ll(1,q^2,\ldots,q^{2(n-1)})$.  Let $\rho = \ll \cup
\mu$ and $\pi = \ll \cup \nu$.  We claim that the summand
$z_\nu^{-1}z_\mu^{-1}z_\ll^{-1}p_\ll(1,q^2,\ldots,q^{2(n-1)})$ is
the coefficient of $p_\mu(\u) p_\nu^\perp(\u)$ when applying
(\ref{eq:hei}) repeatedly to
$z_\pi^{-1}z_\rho^{-1}p_\pi^\perp(\u)p_\rho(\u)$. This is a
straightforward computation, counting the number of ways of
picking parts from $\rho$ and $\pi$ to make the partition $\ll$.

Thus (\ref{eq:hei}) implies Corollary \ref{cor:hhperp}, and since
both the homogeneous and power sum symmetric functions generate
the algebra of symmetric functions, the Corollary follows.


\end{proof}

In the context of the Fock space of $\uqsln$, the operators $B_k$
are known as the \emph{$q$-boson operators} (see \cite{KMS}).  In
\cite{Lam}, we showed that Corollary \ref{cor:KMS} implies Pieri
and Cauchy formulae for the ribbon functions $\g_\ll(X;q)$.

Define the operators $\h{i}{j} = (u_id_i)^j - (d_iu_i)^j$ for $i
,j\in \zz$ and $j \geq 1$. The operators $\h{i}{j}$ act as
follows:
\[
\h{i}{j}: \ll \mapsto \begin{cases} -q^{2j\cdot\spin(\mu/\ll)} \ll
& \mbox{if $\ll$ has an $i$-addable ribbon $\mu/\ll$,} \\
q^{2j\cdot\spin(\ll/\nu)} \ll & \mbox{if $\ll$ has an
$i$-removable ribbon $\ll/\nu$,} \\
    0 &  \mbox{otherwise.} \end{cases}
\]

Since $\h{i}{j}$ acts diagonally on $\f$ in the natural basis,
they all commute with each other.

We will need the following proposition, due to van Leeuwen
\cite{vL}.
\begin{prop}
\label{prop:vl} Let $\ll$ be a partition and suppose ribbons $R_i$
and $R_j$ can be added or removed on diagonals $i < j$ such that
no ribbons can be added or removed on a diagonal $d \in (i,j)$.
Then one of the following holds:
\begin{enumerate}
\item
Both $R_i$ and $R_j$ can be added and $\spin(R_j) = \spin(R_i) -
1$.
\item
Both $R_i$ and $R_j$ can be removed and $\spin(R_j) = \spin(R_i) +
1$.
\item
One of $R_i$ and $R_j$ can be added and the other can be removed
and $\spin(R_j) = \spin(R_i)$.
\end{enumerate}
\end{prop}
For example, writing down, from bottom left to top right, the
spins of the ribbons that can be added and removed from the
partition in Figure \ref{fig:ribbontableau} gives
$2,1,-1,1,-1,1,0$ where positive numbers denote addable ribbons
and negative numbers denote removable ribbons.

\begin{lemma}\label{lem:haction}
Let $\ll$ be a partition, $i \in \zz$ and $j \geq 1$ be an
integer.  Suppose that $\ll$ has an $i$-addable ribbon with spin
$s$.  Then
\[
\brac{\sum_{k={i+1}}^{k=\infty}\h{k}{j}}\ll = -(1 + q^{2j} +
\cdots + q^{2(s-1)j})\ll.
\]
Suppose that a $\ll$ has an $i$-removable ribbon with spin $s$..
Then
\[
\brac{\sum_{k={i+1}}^{k=\infty}\h{k}{j}}\ll = -(1 + q^{2j} +
\cdots + q^{2sj})\ll.
\]
Also, if $i$ is sufficiently small so that no ribbons can be added
to $\ll$ before the $i$-th diagonal, then
\[
\brac{\sum_{k={i}}^{k=\infty}\h{k}{j}}\ll = -(1 + q^{2j} + \cdots
+ q^{2(n-1)j})\ll.
\]
\end{lemma}

\begin{proof}
This follows from Proposition \ref{prop:vl} and the fact that the
furthest ribbon to the right (respectively, left) that can be
added to any partition always has spin $0$ (respectively, $n-1$).
\end{proof}

\begin{lemma}
\label{lem:uhcommute} Let $i,j \in \zz$ and $j \geq 1$.  Then
\[
u_i \brac{\sum_{k={i+1}}^{k=\infty}\h{k}{j}} =
\brac{\sum_{k=i}^{k=\infty}\h{k}{j}}u_i.
\]
Similarly we have,
\[
d_i \brac{\sum_{k={i}}^{k=\infty}\h{k}{j}} =
\brac{\sum_{k=i+1}^{k=\infty}\h{k}{j}}d_i.
\]
\end{lemma}

\begin{proof}
We consider applying the first statement to a partition $\ll$. The
expression vanishes unless $\ll$ has an $i$-addable ribbon
$\mu/\ll$, which we assume has spin $s$.  Then we can compute both
sides using Lemma \ref{lem:haction}.  The second statement follows
similarly.
\end{proof}

We may rewrite equation (\ref{eq:cauchy}) as
\[
U(x)D(y) = D(y)\prod_{i=0}^{n-1} (1-q^ixy)U(x).
\]
Now
\[
\prod_{i=0}^{n-1} (1-q^{2i}xy) = \sum_{i=0}^{n-1} (-1)^i
e_i(1,q^2,\ldots,q^{2(n-1)}) (xy)^i
\]
and $\sum_{i=-\infty}^{i=\infty} \h{i}{j}$ acts as the scalar
$-p_j(1,q^2,\ldots,q^{2(n-1)})= -(1+q^{2j}+\cdots+q^{2(n-1)j})$ by
Lemma \ref{lem:haction}.  Let us use the notation $p_j(h_i) =
-\sum_{k=i}^{\infty} \h{k}{j}$, for $i \in \Z \cup \set{\infty}$.
We also write
\[
e_n(h_i) =\sum_{\rho \vdash n} \epsilon_\rho z_\rho^{-1}
p_{\rho}(h_{i})
\]
where $p_{\rho}(h_i) = p_{\rho_1}(h_i)p_{\rho_2}(h_i)\cdots$ and
$\epsilon_\rho = (-1)^{|\rho|-l(\rho)}$ and $z_\rho$ is as defined
earlier.  As scalar operators on $\f$ we have,
\begin{align*}
\prod_{j=0}^{n-1} (1-q^{2j}xy) &= \sum_{j=0}^{n} (-1)^j
e_j(h_{-\infty})(xy)^j = \sum_{j=0}^{\infty} (-1)^j
e_j(h_{-\infty})(xy)^j.
\end{align*}

\begin{lemma}
\label{lem:cauchycom} Let $i \in \zz$.  Then
\[
(1+yd_i) \left(\sum_{k=0}^{\infty} (-1)^k e_k(h_i)(xy)^k \right)
(1+xu_i) = (1+xu_i) \left( \sum_{k=0}^{\infty} (-1)^k e_k(h_{i+1})
(xy^k)\right) (1+yd_i).
\]
\end{lemma}

\begin{proof}
First we consider the coefficient of $x^{k+1}y^k$.  We need to
show that $e_k(h_i)u_i = u_i e_k(h_{i+1})$ which just follows
immediately from the definition of $e_k(h_i)$ and Lemma
\ref{lem:uhcommute}.  The equality for the coefficient of
$x^ky^{k+1}$ follows in a similar manner.

Now consider the coefficient of $(xy)^k$.  We need to show that
\[
d_ie_{k-1}(h_i)u_i - e_k(h_i) - u_ie_{k-1}(h_{i+1})d_i +
e_{k}(h_{i+1}) = 0.
\]
By Lemma \ref{lem:uhcommute}, this is equivalent to
\begin{equation}
\label{eq:van} \brac{(u_id_i - h_i)e_{k-1}(h_{i+1}) - e_k(h_i) -
u_id_ie_{k-1}(h_i) +e_k(h_{i+1})}\cdot \ll = 0
\end{equation}
for every partition $\ll$.  We now split into three cases
depending on $\ll$.
\begin{enumerate}
\item
Suppose that that $\ll$ has a $i$-addable ribbon with spin $s$.
Then $d_i \cdot \ll = 0$ so (\ref{eq:van}) reduces to
$\brac{e_k(h_{i+1}) +q^{2s} e_{k-1}(h_{i+1}) -e_k(h_i)}\cdot \ll =
0$. Using Lemma \ref{lem:haction}, this becomes
\[
e_k(1,q^2,\ldots,q^{2(s-1)}) + q^{2s} e_{k-1}
(1,q^2,\ldots,q^{2(s-1)}) = e_k(1,q^2,\ldots,q^{2s})
\]
which is an easy symmetric function identity.
\item
Suppose that that $\ll$ has a $i$-removable ribbon with spin $s$.
Then $u_id_i \cdot \ll = h_i \cdot \ll$ so we are reduced to
showing
\[
e_k(h_{i+1}) - e_k(h_i) - q^{2s} e_{k-1}(h_i) = 0
\]
which by Lemma \ref{lem:haction} becomes the same symmetric
function identity as above.
\item
Suppose that $\ll$ has neither a $i$-addable or $i$-removable
ribbon.  Then (\ref{eq:van}) becomes $e_k(h_i)\cdot \ll =
e_k(h_{i+1})\cdot\ll$ so the result is immediate.
\end{enumerate}

\end{proof}

Theorem \ref{thm:cauchy} now follows easily.

\begin{proof}[Proof of Theorem \ref{thm:cauchy}.]
We need to show that  $U(x)D(y) \cdot \ll = D(y)\prod_{i=0}^{n-1}
(1-q^{2i}xy)U(x) \cdot \ll$ for each partition $\ll$.  But if we
focus on a fixed coefficient $x^ay^b$, then for some integers $s <
t$ depending on $\ll$, $a$ and $b$, we may assume that $u_i = d_i
= 0$ for $i < s $ and $i >t$ for all our computations.  Thus it
suffices to show that
\begin{eqnarray*}
&(1 + xu_s)\cdots (1+ xu_t) (1+yd_t) \cdots (1+yd_s)  =\\
&(1+yd_t) \cdots (1+yd_s) \brac{\sum_{k=0}^{\infty} (-1)^k
e_k(h_s)(xy)^k} (1 + xu_s)\cdots (1+ xu_t).\end{eqnarray*} Since
by Lemma \ref{lem:udcommute}, $(1+xu_a)$ and $(1+yd_b)$ commute
unless $a=b$, this follows by applying Lemma \ref{lem:cauchycom}
repeatedly.
\end{proof}

\section{Ribbon functions and $q$-Littlewood Richardson
Coefficients} \label{sec:ribbonfunctions}

In this section we connect the non-commutative Schur functions
$s_\ll(\u)$ with the $q$-Littlewood Richardson coefficients of
Lascoux, Leclerc and Thibon \cite{LLT,LT}.  The set up here is
actually a special case of work of Fomin and Greene \cite{FG},
though not all of their assumptions and results apply in our context.

Let $x_1,x_2, \ldots$ be commutative variables.  Consider the Cauchy product
in the commutative variables $\set{x_i}$ and non-commutative variables
$\set{u_i}$ (not to be
confused with the Cauchy identity in Section \ref{sec:cauchy}):
\[
\Omega(x,\u) = \prod_{j=1}^{\infty}
\left(\prod_{i=\infty}^{\infty} (1+x_j u_i) \right)
\]
where the product is multiplied so that terms with smaller $j$ are
to the right and for the same $j$, terms with smaller $i$ are to
the right.  We have
\begin{align*}
\Omega(x,\u) &= \prod_{j=1}^{\infty} \left( \sum_k x_j^k h_k(\u)
\right) = \sum_\ll s_\ll(X) s_\ll(\u) = \sum_\ll m_\ll(x)
h_\ll(\u)
\end{align*}
where we have used Theorem \ref{thm:hcommute} and the classical
Cauchy identity for symmetric functions.

Since each $h_k(\u)$ adds a horizontal ribbon strip, we see that
\[
\g_{\ll/\mu}(X;q) = \sum_\alpha x^\alpha \br{h_\alpha(\u) \cdot
\mu,\ll} = \sum_\nu m_\nu(x) \br{h_\nu(\u) \cdot \mu,\ll} =
\br{\Omega(x,\u)\cdot \mu,\ll}
\]
where the sum is over all compositions $\alpha$ or all partitions
$\nu$.  In the language of Fomin and Greene \cite{FG}, these
functions were denoted $F_{g/h}$.  The following corollary is
immediate.
\begin{corollary}
The power series $\g_{\ll/\mu}(X;q)$ are symmetric functions in
the (commuting) variables $\set{x_1,x_2,\ldots}$ with coeffcients
in $K$.
\end{corollary}

Letting $\br{.,.}_X$ be the (Hall) inner product in the $X$
variables.  We can write
\begin{align*}
c^\nu_{\ll/\mu}(q) &= \br{\g_{\ll/\mu}(X;q), s_\nu(X)}_X \\
&=  \br{ \br{\Omega(x,\u)\cdot\mu,\ll}, s_\nu(X)}_X \\
&= \br{ \br{\sum_\rho s_\rho(X) s_\rho(\u)\cdot\mu,\ll},s_\nu(X)}_X\\
&= \br{s_\nu(\u)\cdot\mu,\ll}
\end{align*}
so that the action of the noncommutative Schur functions in ribbon
Schur operators $s_\nu(\u)$ calculate the skew $q$-Littlewood
Richardson coefficients.

\begin{prop}
\label{prop:nonneg} The noncommutative Schur function $s_\nu(\u)$
can be written as a non-negative sum of monomials if and only if
the skew $q$-Littlewood Richardson coefficients
$c^\nu_{\ll/\mu}(q) \in {\mathbb N}[q]$ are non-negative
polynomials for all skew shapes $\ll/\mu$.
\end{prop}
\begin{proof}
The only if direction is trivial since $\br{\underline{u}\cdot
\mu,\ll}$ is always a non-negative polynomial in $q$.  Suppose
$c^\nu_{\ll/\mu}(q)$ are non-negative polynomials for all skew
shapes $\ll/\mu$. Write $s_\nu(\u)$ as an alternating sum of
monomials.  Let $\underline{u}$ and $\underline{v}$ be two
monomials occurring in $s_\nu(\u)$.  If $\underline{u}\cdot\mu =
\underline{v} \cdot \mu \neq 0$ for some partition $\mu$ then by
Lemma \ref{lem:commute}, $\underline{u} = \underline{v}$. Now
collect all terms $\underline{v}$ in $s_\nu(\u)$ such that
$\underline{u}\cdot \mu = q^{a(\underline{v})}\underline{v} \cdot
\mu$. Collecting all monomials $\underline{v}$ with a fixed
$a(\underline{v})$, we see that we must be able to cancel out any
negative terms since all coefficients in $c^\nu_{\ll/\mu}(q) =
\br{s_\nu(\u)\cdot\mu,\ll}$ are non-negative.
\end{proof}

To our knowledge, a combinatorial proof of the non-negativity of
the skew $q$-Littlewood Richardson coefficients is only known for
the case $n=2$ via the Yamanouchi domino tableaux of Carr\'{e} and
Leclerc \cite{CL}.  When $\mu = \emptyset$, Leclerc and Thibon
\cite{LT} have shown, using results of \cite{VV} that
$c^\nu_{\ll}(q)$ is equal to a parabolic Kazhdan-Lusztig
polynomial.  Geometric results of \cite{KT} then imply that
$c^\nu_{\ll}(q)$ are non-negative polynomials in $q$.

\section{Non-commutative Schur functions in ribbon Schur operators}
\label{sec:schur} Based on Proposition \ref{prop:nonneg}, we
suggest the following problem, which is the main problem of the
paper.

\begin{conj}
\label{conj:main}  We have:
\begin{enumerate}
\item The non-commutative
Schur functions $s_\ll(\u)$ can be written as a non-negative sum
of monomials.
\item \label{it:strong} A canonical such expression for $s_\ll(\u)$ can be
given by picking some monomials occurring in $h_\ll(\u)$.  That
is, $s_\ll(\u)$ is a positive sum of monomials $\underline{u} =
u_{i_k}u_{i_{k-1}} \cdots u_{i_1}$ where $i_{\ll_1} > i_{\ll_1 -1}
> \cdots > i_1$ and $i_{\ll_1 + \ll_2} > i_{\ll_1+ \ll_2 -1} >
\cdots > i_{\ll_1+1}$ and so on, where $k = |\ll|$.
\end{enumerate}
\end{conj}

By the usual Littlewood-Richardson rule, Conjecture
\ref{conj:main} also implies that the skew Schur functions
$s_{\ll/\mu}(\u)$ can be written as a non-negative sum of
monomials.  If Conjecture \ref{conj:main} is true, we propose the
following definition.

\begin{conjdef}
Let $\ll/\mu$ be a skew shape and $\nu$ a partition so that
$n|\nu| = |\ll/\mu|$.  Let $\underline{u} = u_{i_k}u_{i_{k-1}}
\cdots u_{i_1}$ be a monomial occurring in (a canonical expression
of) $s_\nu(\u)$ so that $\underline{u}.\mu = q^a \ll$ for some
integer $a$.  If Conjecture \ref{conj:main}.(\ref{it:strong})
holds, then the action of $\underline{u}$ on $\mu$ naturally
corresponds to a ribbon tableau of shape $\ll/\mu$ and weight
$\nu$.  We call such a tableau a \emph{Yamanouchi ribbon tableau}.
\end{conjdef}

In their work, Fomin and Greene \cite{FG} show that if the $u_i$
satisfy certain relations then $s_\ll(\u)$ can be written in terms
of the reading words of semistandard tableaux of shape $\ll$.
These relations do not hold, however, for our $u_i$.  We shall see
that a similar description holds for our $s_\ll(\u)$ when $\ll$ is
a hook shape, but appears to fail for other shapes.

The following Theorem holds for any variables $u_i$ satisfying
Theorem \ref{thm:hcommute}.  The commutation relations of
Proposition \ref{prop:commute} are not needed at all.  Let $T$ be
a tableau (not necessarily semistandard).  The \emph{reading word}
$\reading(T)$ is obtained by reading beginning in the top row from
right to left and then going downwards.  If $w = w_1w_2\cdots w_k$
is a word, then we set $u_w = u_{w_1}u_{w_2}\cdots u_{w_k}$.

\begin{thm}
\label{thm:hook} Let $\ll = (a,1^b)$ be a hook shape.  Then
\begin{equation}
\label{eq:hook} s_\ll(\u) = \sum_T u_{\reading(T)}
\end{equation}
where the summation is over all semistandard tableaux $T$ of shape
$\ll$.  For our purposes, the semistandard tableaux can be filled
with any integers not just positive ones, and the row and columns
both satisfy strict inequalities (otherwise $u_{\reading(T)} =
0$).
\end{thm}

\begin{proof}
The theorem is true by definition when $b=0$.  We proceed by
induction on $b$, supposing that the theorem holds for partitions
of the form $(a,1^{b-1})$.  Let $\ll = (a,1^b)$.  The Jacobi-Trudi
formula gives
\[
s_\ll(\u) = \det \left( \begin{array}{ccccc} h_a(\u) & h_{a+1}(\u)
&\cdots & h_{a+b-1}(\u) & h_{a+b}(\u) \\ 1 & h_1(\u) & \cdots &h_{b-1}(\u) & h_{b}(\u) \\
0 & 1 & \cdots & h_{b-2}(\u) & h_{b-1}(\u) \\
\vdots & \vdots & \vdots & \vdots &\vdots \\
0 & 0 &\cdots & 1 & h_1(\u)
\end{array} \right).
\]

Expanding the determinant beginning from the bottom row we obtain
\[
s_\ll(\u) = \sum_{j=1}^{j=b} \left((-1)^{j+1} s_{(a,1^{b-j})}(\u)
h_j(\u)\right) + (-1)^{b} h_{a+b}(\u).
\]
Using the inductive hypothesis, $s_{(a,1^{b-j})}(\u) h_j(\u)$ is
the sum over all the monomials $\underline{u} =
u_{i_1}u_{i_2}\cdots u_{i_{a+b}}$ satisfying $i_1 > i_2
> \cdots > i_a < i_{a+1} < \cdots < i_{a+b-j}$ and $i_{a+b-j+1} >
i_{a+b-j+2} > \cdots > i_{a+b}$.  Let $A_j$ be the sum of those
monomials also satisfying $i_{a+b-j} < i_{a+b-j+1}$ and $B_j$ be
the sum of those such that $i_{a+b-j} > i_{a+b-j+1}$ so that
$s_{(a,1^{b-j})}(\u) h_j(\u) = A_j + B_j$.  Observe that $B_j =
A_{j+1}$ for $j \neq b$ and that $B_{b} = h_{a+b}(\u)$. Cancelling
these terms we obtain $s_\ll(\u) = A_1$ which completes the
inductive step and the proof.
\end{proof}

We speculate that when $n > 2$, the formula (\ref{eq:hook}) of
Theorem \ref{thm:hook} holds if and only if $\ll$ is a hook shape.
We now describe $s_\ll(\u)$ for shapes $\ll$ of the form $(s,2)$.
We will only need the following definition for these shapes, but
we make the general definition in the hope it may be useful for
other shapes.

\begin{definition}
Let $T: \set{(x,y) \in \ll} \rightarrow \zz$ be a filling of the
squares of $\ll$, where $x$ is the row index, $y$ is the column
index and the numbering for $x$ and $y$ begins at 1. Then $T$ is a
\emph{$n$-commuting} tableau if the following conditions hold:
\begin{enumerate}
\item
All rows are increasing, that is, $T(x,y) < T(x,y+1)$ for
$(x,y),(x,y+1) \in \ll$.
\item
If $(x,y), (x+1,y) \in \ll$ and $y > 1$ then $T(x,y) \leq
T(x+1,y)$. Also if $y = 1$ and $(x+1,2) \notin \ll$ then $T(x,1)
\leq T(x+1,1)$.
\item
Suppose the two-by-two square $(x,1),(x,2),(x+1,1),(x+1,2)$ lies
in $\ll$ for some $y$.  Then if $T(x,1) > T(x+1,1)$ we must have
$T(x+1,2)-T(x+1,1) \leq n$.  Otherwise $T(x,1) < T(x+1,1)$ and
either we have $T(x,2) \leq T(x+1,1)$ or we have both $T(x,2) >
T(x+1,1)$ and $T(x+1,2) - T(x,1) \leq n$.
\end{enumerate}
\end{definition}

We give some examples of commuting and non-commuting tableaux of
shape $(2,2)$ in Figure \ref{fig:nice}.

\begin{figure}
 \pspicture(20,0)(100,120)
\abox(1,1){2}\abox(2,1){3}\abox(1,2){0}\abox(2,2){1}
\abox(4,1){1}\abox(5,1){4}\abox(4,2){0}\abox(5,2){2}
\abox(7,1){0}\abox(8,1){3}\abox(7,2){1}\abox(8,2){2}

\abox(1,5){0}\abox(2,5){1}\abox(1,6){2}\abox(2,6){3}
\abox(4,5){1}\abox(5,5){3}\abox(4,6){0}\abox(5,6){2}
\abox(7,5){0}\abox(8,5){4}\abox(7,6){1}\abox(8,6){2}

\endpspicture
\caption{The bottom row contains some 3-commuting tableaux and the
top row contains some tableaux which are not 3-commuting.}
\label{fig:nice}
\end{figure}

Let us now note that the operators $u_i$ satisfy the following
Knuth-like relations, using Proposition \ref{prop:commute}. Let
$i,j,k$ satisfy either $i < j < k$ or $i
> j > k$, then depending on whether $u_i$ and $u_k$ commute, we
have
\begin{enumerate}
\item
\label{it:knuth1} Either $u_i u_k u_j = u_k u_i u_j$ or $u_i u_k
u_j = u_j u_i u_k$.
\item
\label{it:knuth2} Either $u_j u_k u_i = u_j u_i u_k$ or $u_j u_k
u_i = u_k u_i u_j$.
\end{enumerate}
Note that the statement of these relations do not explicitly
depend on $n$.

\begin{thm}
\label{thm:a2} Let $\ll =(s,2)$.  Then
\[
s_\ll(\u) = \sum_T u_\reading(T)
\]
where the summation is over all $n$-commuting tableaux $T$ of
shape $\ll$.
\end{thm}

\begin{proof}
By definition $s_\ll(\u) = h_s(\u)h_2(\u) - h_{s+1}(\u)h_1(\u)$.
We have that $h_{s+1}(\u)h_1(\u)$ is equal to the sum over
monomials $\underline{u} = u_{i_1}\cdots u_{i_{s+1}}u_{i_{s+2}}$
such that $i_1
> i_2 > \cdots > i_{s+1}$.  We will show how to use the Knuth-like
relations (\ref{it:knuth1}) and (\ref{it:knuth2}) to transform
each such monomial into one that occurs in $h_s(\u)h_2(\u)$, in an
injective fashion.

Let $a = i_{s-1}$, $b = i_{s}$, $c = i_{s+1}$ and $d = i_{s+2}$ so
that $a > b > c$.  We may assume that $a,b,c,d$ are all different
for otherwise $\underline{u} = 0$ by Prop \ref{prop:commute}.  If
$c > d$ then $\underline{u}$ is already a monomial occurring in
$h_s(\u)h_2(\u)$.  So suppose $c < d$. Now if $d>  b > c$ we apply
the Knuth-like relations to transform $\underline{v} = u_a u_b u_c
u_d$ to either $\underline{v}' = u_a u_b u_d u_c$ or
$\underline{v}' = u_a u_c u_d u_b$, where we always pick the
former if $u_d$ and $u_c$ commute.  If $b
> d > c$ we may transform $\underline{v} = u_a u_b u_c u_d$ to either $\underline{v}'=u_a u_c u_b
u_d$ or $\underline{v}'= u_a u_d u_b u_c$ again picking the former
if $u_b$ and $u_c$ commute. In all cases the resulting monomial
occurs in $h_s(\u) h_2(\u)$ and the map
$\underline{v}\mapsto\underline{v}'$ is injective if the
information of whether $u_c$ commutes with both $u_b$ and $u_d$ is
fixed. Writing $\set{j_0<j_1<j_2<j_3}$ for $\set{a,b,c,d}$ to
indicate the relative order we may tabulate the possible resulting
monomials $\underline{v}'$ as reading words of the following
``tableaux'':
\begin{align*}
a > b > c > d:\;& \tbt{2}{3}{0}{1} \mbox{($j_0< j_1<j_2< j_3 \in \Z$)};&\\
a > b > d > c:\; & \tbt{0}{3}{1}{2} \mbox{($u_{j_0}$ and $u_{j_2}$ commute)} & \tbt{1}{3}{0}{2} \mbox{($u_{j_0}$ and $u_{j_2}$ do not commute)} ;&\\
a > d > b > c:\;& \tbt{1}{3}{0}{2}  \mbox{($u_{j_0}$ and $u_{j_2}$ commute)}& \tbt{0}{3}{1}{2} \mbox{($u_{j_0}$ and $u_{j_2}$ do not commute)};&\\
d > a > b > c:\; & \tbt{1}{2}{0}{3} \mbox{($u_{j_0}$ and $u_{j_3}$
commute)}& \tbt{0}{2}{1}{3} \mbox{($u_{j_0}$ and $u_{j_3}$ do not
commute)}.&
\end{align*}
Finally, cancelling these monomials from $h_s(\u) h_2(\u)$, we see
that the monomials in $s_\ll(\u)$ are of the form $u_{i_1}\cdots
u_{i_{s-3}}u_{i_{s-2}}u_{\reading(T)}$ for a tableau $T$ of the
following form, (satisfying $i_{s-2} > t$ where $t$ is the value
in the top right hand corner):
\begin{align*}
\tbt{0}{1}{2}{3} & \ \; \mbox{with no extra conditions,} \\
\tbt{1}{2}{0}{3} &\  \; \mbox{if $u_{j_3}$ and $u_{j_0}$ do not
commute,} \\
\tbt{0}{2}{1}{3} &\  \; \mbox{if $u_{j_3}$ and $u_{j_0}$ commute,}
\end{align*}
where in the first case $j_0 < j_1 \leq j_2 < j_3$ and in the
remaining cases $j_0 < j_1 < j_2 < j_3$.  We may allow more
equalities to be weak, but the additional monomials
$u_{\reading(T)}$ that we obtain turn out all to be 0.  One
immediately checks that these monomials are the reading words of
$n$-commuting tableaux of shape $(s,2)$.
\end{proof}

By Theorem \ref{thm:hook}, $e_k(\u) = \sum_{i_1 < i_2 < \cdots <
i_k} u_{i_1} \cdots u_{i_k}$.  Thus using the dual Jacobi-Trudi
formula, we can get a description for $s_{\ll'}(\u)$ whenever we
have one for $s_{\ll}(\u)$ by reversing the order `$<$' on $\zz$.
This for example leads to a combinatorial interpretation for
$s_\ll(\u)$ of the form $\ll = (2,2,1^a)$ which we will not write
out explicitly.

One can also obtain different combinatorial interpretations for
$s_\ll(\u)$ by for example changing $h_a(\u)h_2(\u)$ to
$h_2(\u)h_a(\u)$ in the proof of Theorem \ref{thm:a2}.  This leads
to the reversed reading order on tableaux which also has the order
`$<$' reversed.  In the case $n=2$, the $n$-commuting tableaux are
nearly the same as usual semistandard tableaux where row and
column inequalities are strict.  The reversed reading order on
order-reversed semistandard tableaux would lead to the same
combinatorial interpretation as Carr\'{e} and Leclerc's Yamanouchi
domino tableaux \cite{CL}.

It seems likely that Theorem \ref{thm:hook} and Theorem
\ref{thm:a2} may be combined to give a description of $s_\ll(\u)$
for $\ll = (a,2,1^b)$ but so far we have been unable to make
progress on the case $\ll = (3,3)$.

We end this section with an example.
\begin{example}
Let $n=3$ and $\ll = (4,4,4)$. Let us calculate the coefficient of
$s_{22}$ in $\g_{\ll}(X;q)$. The shape $\ll$ has ribbons on the
diagonals $\set{0,1,2,3}$, so we are concerned with $3$-commuting
tableaux of shape $(2,2)$ filled with the numbers $\set{0,1,2,3}$.
There are only two such tableaux $S$ and $T$ given in Figure
\ref{fig:nice0123}.  It is easy to see that
$u_{\reading(T)}.\emptyset = 0$ so the coefficient of $s_{22}$ in
$\g_{\ll}(X;q)$ is given by the spin of the ribbon tableau
corresponding to $u_{\reading(S)}.\emptyset = u_2 u_1 u_3
u_0.\emptyset$. This tableau has spin $4$, so that
$c^{(2,2)}_{(4,4,4)}(q) = q^4$ (see Figure \ref{fig:spin4}).

\begin{figure}
 \pspicture(20,0)(100,60)

\rput(-5,30){$T = $}
\rput(95,30){$S = $}
\abox(1,1){2}\abox(2,1){3}\abox(1,2){0}\abox(2,2){1}
\abox(7,1){0}\abox(8,1){3}\abox(7,2){1}\abox(8,2){2}
\endpspicture
\caption{The 3-commuting tableaux with shape $(2,2)$ and squares
filled with $\set{0,1,2,3}$.} \label{fig:nice0123}
\end{figure}

\begin{figure}
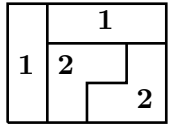


\pspicture(0,10)(100,50)
\psline(0,0)(0,45)(60,45)(60,0)(0,0)
\psline(15,0)(15,45)
\psline(15,30)(60,30)
\psline(30,0)(30,15)(45,15)(45,30)
\rput(7,22){\bf 1}
\rput(22,22){\bf 2}
\rput(37,39){\bf 1}
\rput(52,9){\bf 2}
\endpspicture
\caption{The Yamanouchi ribbon tableau corresponding to
$c^{(2,2)}_{(4,4,4)}(q) = q^4$.} \label{fig:spin4}
\end{figure}

We may check directly that in fact
\[
\g_{\ll}(X;q) = q^2 s_{211} +q^4(s_{31}+s_{22}) + q^6 s_{31} + q^8 s_4.
\]
\end{example}

\section{Final remarks}
\label{sec:final}
\subsection{Maps involving $\U_n$}
The algebra
$\U_n$ has many automorphisms.  The map sending $u_i \mapsto
u_{i+1}$ is an algebra isomorphism of $\U_n$ and the map sending
$u_i \mapsto u_{-i}$ a semi-linear algebra involution.

\begin{prop}
There are (many) commuting injections
\[
\U_1 \hookrightarrow \U_2 \hookrightarrow \U_3 \cdots
\]
and surjections
\[
\cdots \U_3 \twoheadrightarrow \U_2 \twoheadrightarrow \U_1.
\]
\end{prop}
\begin{proof}
The injection $\U_{n-1} \hookrightarrow \U_n$ can be given by
sending $u_{k(n-1)+i}$ to $u_{kn+i}$ where $i \in
\set{0,1,\ldots,n-2}$. The surjection $\U_{n} \twoheadrightarrow
\U_{n-1}$ is given by sending $u_{kn+i}$ to $u_{k(n-1)+i}$ for $i
\in \set{0,1,\ldots,n-2}$ and sending $u_{kn+n-1}$ to 0 for all $k
\in \zz$.
\end{proof}

\subsection{The algebra $\U_\infty$} Picking compatible injections as above, the inductive limit $\U_\infty$ of
the algebras $\U_n$ has a countable set of generators $u_{i,j}$
(the image of $u_{i + jn}^{(n)}$) where $i \in \nn$ and $j \in
\zz$. The generators are partially ordered by the relation $(i,j)
< (k,l)$ if either $j < l$, or $j = l$ and $i < k$.  The
generators satisfy the following relations
\begin{align*}
u_{i,j}u_{k,l}&=u_{k,l}u_{i,j} &\mbox{if $|j - l| \geq 2$ or if $j = l \pm 1$ and $i \neq k$,} \\
u_{i,j}^2 & = 0& \mbox{for any $i,j \in \Z$,}\\
u_{i,j}u_{i,j \pm 1}u_{i,j} &= 0 & \mbox{for any $i,j \in \Z$,}\\
u_{i,j} u_{k,j} &= q^2 u_{k,j} u_{i,j} &\mbox{for $i,j,k \in \Z$  satisfying $i >k$,} \\
u_{i,j + 1} u_{k,j} &= q^2 u_{k,j} u_{i,j + 1} &\mbox{for $i,j,k
\in \Z$ satisfying $i < k$.}
\end{align*}

Many of the results of this paper can be phrased in terms of
$\U_\infty$.  For example, one can define $\infty$-commuting
tableaux which are maps $T:\ll \rightarrow \nn \times \zz$.

\subsection{Another description of $\U_n$} There is an alternative way
of looking at the algebras $\U_n \subset \mathrm{End}(\f)$. Let $u
= u_0$ be the operator adding a $n$-ribbon on the $0$-th diagonal.
We may view a partition $\ll$ in terms of its
$\set{0,1}$-\emph{edge sequence} $\set{p_i(\ll)}_{i =
-\infty}^{\infty}$ (see \cite{EC2} for example). Let $t$ be the
operator which shifts a bit sequence
$\set{p_i}_{i=-\infty}^{\infty}$ one-step to the right, so that
$(t\cdot p)_i = p_{i-1}$.  Note that $t$ does not send a partition
to a partition so we need to consider it as a linear operator on a
larger space (spanned by doubly-infinite bit sequences, for
example).  It is clear that $u_i = t^{i}ut^{-i}$ and we may
consider $\U_n$ as sitting inside an enlarged algebra
$\tilde{\U_n} = K[u,t,t^{-1}]$.

\subsection{Relationship with affine Hecke algebras and canonical bases}  The operators $p_k(\u)=B_i$ of
Corollary \ref{cor:KMS} were originally described in \cite{KMS} as
power sums in certain elements $y_i^{-1}$ which lie within the
thermodynamic limit of the center of the affine Hecke algebra of
type $A$.  In fact, we have $s_\ll(\u) =
s_\ll(y_1^{-1},y_2^{-1},\ldots)$ as operators acting on $\f$.  It
is not clear what the relationship between these $y_i^{-1}$ and
the $u_i$ are. It is also not always true that the action of a
monomial in the $y_i^{-1}$ on a partition $\ll$ gives a power of
$q$ times some other partition (or 0), which is the case for
monomials in $u_i$. Separately, Leclerc and Thibon \cite{LT} have
shown that $s_\ll(\u).\emptyset = G_{n\ll}^-$ a member of the
lower global crystal basis of $\f$. It would be interesting to
understand this in terms of the $u_i$.

\subsection{Connection with work of van Leeuwen and Fomin} Van Leeuwen \cite{vL} has given a spin-preserving
Robinson-Schensted-Knuth correspondence for ribbon tableaux.  The
calculations in Section \ref{sec:cauchy} are essentially an
algebraic version of van Leeuwen's correspondence, in a manner
similar to the construction of the generalised Schensted
correspondences in \cite{Fom}.

\end{document}